\documentclass[12pt]{article}
\newtheorem{theoreme}{Th\'eor\`eme}[section]
\newtheorem{lemma}[theoreme]{Lemme}
\newtheorem{corollary}[theoreme]{Corollaire}
\newtheorem{definition}[theoreme]{D\'efinition}

\newcommand{\Oe}{\emptyset}

\begin{document}

\title{Sur l'existence d'une cat\'egorie ayant une matrice strictement positive donn\'ee}

\author{Samer Allouch\\Universit\'e de Nice-Sophia Antipolis}

\date{Mai, 2008}

\maketitle

\section{Introduction}
Depuis quelques dizaines d'ann\'ees, certains math\'ematiciens se
dirigent vers les domaines de la g\'eom\'etrie alg\'ebrique
reli\'es sp\'ecialement aux cat\'egories finies qui prennent ainsi
une importance croissante dans les math\'ematiques pures et
appliqu\'ees, voir \cite{CuntzHeckenberger} par exemple. 
Le but de mon recherche est d'\'etudier la classification des cat\'egories
finies au moyen de la correspondance entre
les cat\'egories finies d'ordre $n$ et les matrices carr\'ees de
taille $n$. Cette correspondance figure dans les papiers r\'ecents
\cite{BergerLeinster} \cite{Leinster}, voir aussi \cite[p. 486]{Kapranov}.
Notre premi\`ere
t\^ache consiste \`a
trouver quelles matrices correspondent effectivement \`a des cat\'egories.

Pour chaque cat\'egorie $A$ qui a $n$ objets
$x_1$,$x_2$,...,$x_n$, la matrice $M$ de taille $n$
associ\'ee \`a $A$ est d\'efinie par $m_{ij}:= | A(x_i,x_j)|$. Dans le
sens contraire on ne peut pas dire que pour chaque matrice $M$ il
y a une cat\'egorie. Dans mon m\'emoire \cite{Allouch} et dans
\cite{BergerLeinster}, un premier exemple  d'une matrice $M$ qui n'a pas de
cat\'egorie est
\begin{displaymath} \mathbf{M} = \left(
\begin{array}{cc}
2 & 2  \\
1 & 1
\end{array} \right) .
\end{displaymath}

On cherche donc \`a savoir quelles sont les matrices qui ont au moins
une cat\'egorie associ\'ee.
On notera cette condition $Cat(M)\neq{\Oe}$ et on dira parfois que {\em $M$ marche}, ou dans le cas contraire que {\em $M$
ne marche pas} si $Cat(M)={\Oe}$. \\
Donc le but de ce rapport est de trouver une cat\'egorie pour
chaque type de matrice si possible, ou sinon de pouver  que
$Cat(M)={\Oe}$. Pour cel\`a il suffit de faire la d\'emonstration
pour $n=2$ et $n=3$ et ensuite  on
g\'en\'eralise sur $n$. Nous ne traitons ici que le cas o\`u les coefficients sont strictement positifs. \\
Pour $n=2$:\\
On a une matrice $M$ d\'efinie par:
\begin{displaymath}
\mathbf{M} = \left( \begin{array}{cc}
a & b  \\
c & d
\end{array} \right)
\end{displaymath}
avec $a,b,c,d\geq 1$.\\
Pour qu'on a  $Cat(M)\neq{\Oe}$ (il existe au moins une
cat\'egorie associ\'ee \`a $M$), on suppose que soit $a>1$ et
$d>1$, soit $a>bc$ ou $d>bc$, soit $a=b=c=d=1$.
\\
Pour n=3:\\
Soit $M$ une matrice d\'efinie par:
\begin{displaymath}
\mathbf{M} = \left( \begin{array}{ccc}
m_{11} & m_{12}& m_{13} \\
m_{21} & m_{22}&m_{23}\\
m_{31} & m_{32}&m_{33}
\end{array}
\right).
\end{displaymath}
avec $m_{ij}\geq 1$ pour tout $i,j\in\{1,2,3\}$. \\
Si $m_{ii}>1$ pour $i=1,2,3$ on trouve une cat\'egorie $A$
associ\'ee \`a $M$ (Th\'eor\`eme \ref{leinsterthm}). Si par
exemple $m_{11}=1$, soient $m_{22}>m_{12}m_{21}$,
$m_{33}>m_{13}m_{31}$ et $m_{23}\geq m_{21} m_{13}$, $m_{32}\geq
m_{31}m_{12}$ on va d\'emontrer $Cat(M)\neq{\Oe}$.
\\
Pour $n>3$ alors soit $M$ d\'efinie par:
\begin{displaymath}
\mathbf{M} = \left( \begin{array}{cccc}
m_{11} & m_{12} & \ldots &m_{1n}\\
m_{21} & m_{22} & \ldots & m_{2n} \\
\vdots & \vdots & \ddots&\vdots\\
m_{n1} & m_{n2} &  \ldots   &m_{nn}
\end{array} \right)
\end{displaymath}
avec $m_{ij}> 0$ pour tout $i,j \in\{1,2,...,n\}$.\\
Premier cas :\\
S'il existe une seule coefficient diagonale $m_{aa}=1$ alors
$Cat(M)\neq{\Oe}$ si et seulement si
$m_{ii}> m_{ai}m_{ia}\forall (i\neq a)$, $m_{ij}\geq m_{ia}m_{aj} \forall (i\neq j)$.\\
Deuxi\`eme cas:\\
S'il existe plus qu'une coefficient diagonale \'egale \`a $1$ alors il y a deux cas (voir la d\'efinition \ref{defreduite}):\\
-si $M$ est r\'eduite alors M ne marche pas.\\
-si $M$ est non r\'eduite, on r\'eduit $M$ en une sous-matrice $N$, et dans ce
cas $M$
marche si et seulement si $N$ marche.\\
Dernier cas: \\
Si $m_{ii}>1 ,\forall i$  alors $Cat(M)\neq{\Oe}$ (Th\'eor\`eme \ref{leinsterthm}).\\
A la fin nous pouvons dire: Si $M$ est une matrice de taille
$n\geq 3$ avec $m_{ij}\geq 1$, alors $Cat(M)\neq \Oe$ si et
seulement si, pour toute sous-matrice $N\subset M$ de taille $3$
on a $Cat(N)\neq \Oe$.\\
Je veux exprimer mon remerciement plus profonde \`a mon directeur
Carlos Simpson pour sa confiance en mes capacit\'es et mon
travail. Merci de consacrer vos heures pr\'ecieuses pour me
guider.
\section {Matrices carr\'ees d'order 2}
On va \'etudier les matrices de taille 2:\\
\begin{displaymath}
\mathbf{M} = \left( \begin{array}{cc}
a & b  \\
c & d
\end{array} \right) .
\end{displaymath}
\begin{theoreme}
Soit $M$ une matrice carr\'ee d'order 2 d\'efinie par:\\
\begin{displaymath}
\mathbf{M} = \left( \begin{array}{cc}
a & b  \\
c & 1
\end{array} \right)
\end{displaymath} avec $a,b,c>1 $, alors
$Cat(M)\neq{\Oe}$  $\Leftrightarrow|M|=(a-bc)>0$ .
\end{theoreme}
Il faut d\'emontrer les deux sens. \\
a) $Cat(M)\neq{\Oe} \Rightarrow |M|=(a-bc)>0$.\\
En effet: soit $A$ une cat\'egorie finie d'order 2 associ\'ee \`a $M$
dont les objets
sont $\{x_1,x_2\}$ et les morphismes  d\'efinis par:\\
$A(x_1,x_1)=\{e_1=id_{x_1},e_2,...,e_a\}$.\\
$A(x_1,x_2)=\{f_1,f_2,...,f_b\}$.\\
$A(x_2,x_1)=\{g_1,g_2,...,g_c\}$.\\
$A(x_2,x_2)=\{id_{x_2}\}$.\\
\\
Les \'equations de la loi de composition sont :\\
1-$e_{i}e_{j}$ pour tout $i,j\in\{1,...,a\}$.\\
2-$e_i^2$pour tout $i\in\{1,...,a\}$.\\
3-$f_j e_i$ pour tout $i\in\{1,...,a\}$ et $j\in\{1,...,b\}$.\\
4-$g{_i}f_j$ pour tout $i\in\{1,...,c\}$ et$j\in\{1,...,b\}$.\\
5-$f_{i}g_j=id_{x_2}$ pour tout $i\in\{1,...,b\}$ et$j\in\{1,...,c\}$.\\
6-$e_{i}g_j$ pour tout $i\in\{1,...,a\}$ et $j\in\{1,...,c\}$.\\
\\
Normalement il y a huit \'equations d'associativit\'e soit: \\
a-$e_i^3$ vraie pour tout $i\in\{1,2,...,a\}$.\\
b-$f_j(e_i^2)=(f_{j}e_i)e_i$ pour tout $i\in\{1,...,a\}$ et$j\in\{1,...,b\}$.\\
c-$(g_{n}f_j)e_i=g_n(f_{j}e_i)$ pour tout $i\in\{1,...,a\}$ ,$j\in\{1,...,b\}$ et $n\in\{1,...,c\}$.\\
d-$(e_{i}g_n)f_j=e_i(g_{n}f_j)$ pour tout $i\in\{1,...,a\}$ ,$j\in\{1,...,b\}$ et $n\in\{1,...,c\}$.\\
e-$(f_ig_n)f_j=f_i(g_{n}f_j)$ pour tout $i,j\in\{1,...,b\}$ , $n\in\{1,...,c\}$.\\
f-$e_i^2g_n=e_i(e_i g_n)$ pour tout $i\in\{1,...,a\}$  et $n\in\{1,...,c\}$.\\
g-$(f_j e_i)g_n=f_j(e_i g_n)$ pour tout $i\in\{1,...,a\}$ ,$j\in\{1,...,b\}$ et $n\in\{1,...,c\}$.\\
h-$(g_n f_j)g_m=g_n(f_j g_m)$ pour tout $j\in\{1,...,b\}$ et $n,m\in\{1,...,c\}$.\\
Voici quelques remarques:\\
\textbf{Rq (1)}:\\
Soit $(f,g)\in\{f_1,f_2,...,f_b\}\times\{g_1,g_2,...,g_c\}$ tel
que $gf=id_{x_1}$ on a $feg=id_{x_2} \Longrightarrow g(feg)=g$
alors $(gf)(eg)=g\Longrightarrow id_{x_1}(eg)=g$ ce qui donne
$eg=g$ avec ($gf=id_{x_1}$).\\
Aussi $(eg)f=gf=id_{x_1}=e(gf)=e$ impossible donc $f,g$ n'existent
pas tel que $gf=id_{x_1}$.\\
\textbf{Rq( 2)}:\\
$e_i^2=e_i$ pour tout $i\in\{1,...,a\}$ \`a condition qu'il existe $(f,g)\in\{f_1,f_2,...,f_b\}\times\{g_1,g_2,...,g_c\}$ tel que $gf=e_i$.\\
En effet: $fg=id_{x_2}\Longrightarrow g(fg)=g id_{x_2}=g$ alors $(gf)g=g \Longrightarrow eg=g$\\
donc pour tout $eg=g$ pour tout $e\in\{e_2,...,e_a\}$ tel que
$gf=e$.\\
On a $(eg)f=e(gf)\Longrightarrow gf=ee=e^2$.\\
\textbf{Rq (3)}:\\
Il n'existent pas $f_1\neq f_2$ et $g$ tel que $gf_1=gf_2=e$. \\
En effet : $(f_1g)f_2=id_{x_2}f_2=f_2$\\
autrement dit:$(f_1g)f_2=f_1(gf_2)=f_1(gf_1)=(f_1g)f_1=f_1$\\
donc $f_1=f_2$ contradiction.\\
\textbf{Rq (4)}:\\
Il n'existent pas $f$ et $g_1\neq g_2$ tel que $g_1f=g_2f=e$ \\
En effet : $(g_1f)g_2=eg_2=g_2=g_1=g_1(fg_2)$ alors $g_1=g_2$  contradiction.\\
\textbf{Rq (5)}:\\
Il n'existent pas $g_1\neq g_2$ et $f_1\neq f_2$ tel que $g_1f_1=g_2f_2=e$ \\
En effet : $(f_1g_2)f_2=id_{x_2}f_2=f_2$\\
autrement dit:$(f_1g_2)f_2=f_1(g_2f_2)=f_1(g_1f_1)=(f_1g_1)f_1=f_1$\\
donc $f_1=f_2$ contradiction.\\
\\
Donc les 4 remarques \textbf{Rq (1)} \textbf{Rq (3)} \textbf{Rq (4)} et \textbf{Rq (5)} donnent $(a-bc)>0$ c.\`a.d $a>bc$.\\
\\
b) $a>bc\Rightarrow Cat(M)\neq{\Oe}$. \\
On commence par une cat\'egorie dont les \'equations de la loi de composition sont :\\
\\
1-$e_{j}^i e_{p}^n=e_j^n (e^2=e)$ pour tout $i,n\in\{1,...,b\}$ et $j,p\in\{1,...,c\} $.\\
\\
2-$f_n e_j^i=f_i$ pour tout $i,n\in\{1,...,b\}$ et $j\in\{1,...,c\}$.\\
\\
3-$g_j f_i=e_j^i$ pour tout $i \in\{1,...,b\}$ et$j\in\{1,...,c\}$.\\
\\
4-$f_{i}g_j=id_{x_2}$ pour tout $i\in\{1,...,b\}$ et$j\in\{1,...,c\}$.\\
\\
5-$e_j^ig_k=g_j$ pour tout $i\in\{1,...,a\}$ et $j\in\{1,...,c\}$.\\
\\
On va verifie les \'equations associatives.\\
$(e_i^j e_p^n) e_r^s=e_i^n e_r^s=e_i^s=e_i^j (e_p^n e_r^s)=e_i^j e_p^s=e_i^s$ vraie.\\
\\
$(f_j e_r^s)e_p^n=f_s e_p^n=f_n=f_j(e_r^s e_p^n)=f_j e_r^n=f_n$ vraie .\\
\\
$ (g_j f_i)e_p^n=e_j^i e_p^n=e_j^n=g_j(f_i e_p^n)=g_j f_n=e_j^n$  vraie.\\
\\
$(e_p^n g_j)f_i=g_p f_i=e_p^i=e_p^n(g_j f_i)e_p^n e_j^i=e_p^i$ vraie .\\
\\
$(f_j g_n)f_p=f_p=f_j(g_n f_p)=f_j e_n^p=f_p$ vraie .\\
\\
$(e_p^n e_r^s)g_j=e_p^s g_j=g_p=e_p^n(e_r^s g_j)=e_p^n g_r=g_p$ vraie .\\
\\
$(f_i e_p^n)g_j=f_n g_j=id_{x_2}=f_i (e_p^n g_j)=f_i g_p=id_{x_2}$ vraie .\\
\\
$(g_j f_i) g_n=e_j^i g_n=g_j=g_j(f_i g_n)=g_j$ vraie .\\
\\
Toutes les \'equations associatives marchent donc $B$ est une
cat\'egorie finie d'order 2 , d'apres la construction des
morphismes de $B$ alors $B$
est associ\'ee \`a $M$ d\'efinie par :\\
\begin{displaymath}
\mathbf{M} = \left( \begin{array}{cc}
bc+1 & b  \\
c   & 1
\end{array} \right) ,
\end{displaymath}
finalement $Cat(M)\neq{\Oe}$.\\
On note $M(u)$ par:
\begin{displaymath}
\mathbf{M(u)} = \left( \begin{array}{cc}
u & b  \\
c   & 1
\end{array} \right)
\end{displaymath}
avec pour commencer $u=bc+1$.\\
En fait $A$ est une cat\'egorie de $M(u)$, soit $A'$ une nouvelle
cat\'egorie dont les objets $Ob(A)=Ob(A')$ et les morphismes
d\'efinis par :\\
$A'(x_1,x_1)=\{id_{x_1},e^1_1,...,e_j^i,...,e_c^b,e_1,...,e_n\}$\\
$A'(x_1,x_2)=\{f_1,f_2,...,f_b\}$\\
$A'(x_2,x_1)=\{g_1,g_2,...,g_c\}$\\
$A'(x_2,x_2)=\{id_{x_2}\}$\\
avec la loi de composition d\'efinie par:\\
$e_ig_j=e_1^1 g_j = g_1$ pour tout $i\in\{1,2..,n\}$ .\\
\\
$f_i e_j=f_i e_c^b = f_b$ pour tout $i\in\{1,2..,n\}$. \\
\\
$e_i e_j=e_1^1 e_c^b=e_1^b$ pour tout $i,j \in\{1,2..,n\}$ .\\
Toutes les \'equations associatives marchent alors $A'$ est
une cat\'egorie associ\'e \`a $M(u+n)$.
\begin{corollary}
Soit $M$ une matrice carr\'e d'order 2 d\'efinie par:
\begin{displaymath}
\mathbf{M} = \left( \begin{array}{cc}
a & b   \\
c & d
\end{array} \right) ,
\end{displaymath}
alors  on a $Cat(M)\neq{\Oe}$ dans les cas suivants :\\
1- a=b=c=d=1 voir \cite{Allouch};\\
2- a=1, $d>bc$;\\
3- d=1, $a>bc$;\\
4- $a>1$ , $d>1$ voir \cite{Leinster}.\\
Sinon $Cat(M)={\Oe}$.
\end{corollary}
Preuve: on pose $a>bc$.\\
Soit $N$ une matrice d\'efinie par:
\begin{displaymath}
\mathbf{N} = \left( \begin{array}{cc}
a & b   \\
c & 1
\end{array} \right)
\end{displaymath}
D'apr\`es ce qui pr\'ec\`ede il existe $A$ une cat\'egorie associ\'ee
\`a $N$ , avec $Ob(A)=\{x_1,x_2\}$ et les morphismes d\'efinis
par:\\
$A(x_1,x_1)=\{id_{x_1},e^1_1,...,e_j^i,...,e_c^b,e_{bc+2}...,e_a\}$.\\
$A(x_1,x_2)=\{f_1,f_2,...,f_b\}$.\\
$A(x_2,x_1)=\{g_1,g_2,...,g_c\}$.\\
$A(x_2,x_2)=\{id_{x_2}\}$.\\
Soit $A'$ une nouvelle cat\'egorie tel que $Ob(A')=Ob(A)$ et les
morphismes d\'efinis par:\\
$A'(x_1,x_1)=\{id_{x_1},e^1_1,...,e_j^i,...,e_c^b,e_{bc+2}...,e_a\}$.\\
$A'(x_1,x_2)=\{f_1,f_2,...,f_b\}$.\\
$A'(x_2,x_1)=\{g_1,g_2,...,g_c\}$.\\
$A'(x_2,x_2)=\{id_{x_2},n_1,...,n_{d-1}\}$.\\
avec $n_i\neq n_j$ pour tout $i,j$ .\\
Les \'equations de la loi de composition qui d\'ependent des $\{n_1,..,n_{d-1}\}$ sont\\
$nn'=n_1$\\
$nf=f$\\
$gn=g$.\\
Toutes les \'equations associatives marchent donc $A'$ est une cat\'egorie
associ\'ee \`a $M$, donc $Cat(M)\neq{\Oe}$.

\section {Matrices triples}
On va \'etudier les matrices de taille 3:
\begin{displaymath}
\mathbf{M} = \left( \begin{array}{ccc}
1 & a &b  \\
c & n &m\\
p & q &r
\end{array} \right)
\end{displaymath}
Supposons $n>1$ et $r>1$, sinon voir le cas des lemmes (\ref{reduire}, \ref{marchepas}). On peut
d\'eduire de l'\'etude pr\'ec\`edant que $n\geq ac+1$ et $r\geq
bp+1$. Les deux matrices suivantes
sont sous-matrices de $M$:
\begin{displaymath}
\left( \begin{array}{cc}
1 & b  \\
p & r
\end{array} \right)
\end{displaymath}
et
\begin{displaymath}
\left( \begin{array}{cc}
1 & a  \\
c & n
\end{array} \right) ,
\end{displaymath}
voir \cite {Allouch} (ici une {\em sous-matrice} correspond par convention
au m\^eme sous-ensemble des colonnes que des lignes). Un lemme facile dit que si $N$ est une sous matrice qui
ne marche pas alors $M$ ne marche pas non plus !\\
En plus $m\geq bc$ et $q\geq ap$ en effet:\\
soit $A$ une cat\'egorie associ\'ee \`a M  dont les objets sont $\{x_1,x_2,x_3\}$ et les morphismes sont donn\'es par:\\
$A(x_1,x_1)=\{1\}$.\\
$A(x_1,x_2)=\{f_1,...,f_a\}$.\\
$A(x_1,x_3)=\{h_1,...,h_b\}$.\\
$A(x_2,x_1)=\{g_1,...,g_c\}$.\\
$A(x_2,x_2)=\{e_0=1,...,e_n\}$.\\
$A(x_2,x_3)=\{k_1,...,k_m\}$.\\
$A(x_3,x_1)=\{L_1,...,L_p\}$.\\
$A(x_3,x_2)=\{M_1,...,M_q\}$.\\
$A(x_3,x_3)=\{N_1,...,N_r\}$.\\
\\
les \'equations de la loi de composition sont:
$$\begin{tabular}{|c|c|c|c|c|c|}
\hline$NN'\in\{ N \}$&$MN\in\{M\}$\\
\hline$LN\in\{L\}$&$gM\in\{L\}$\\
\hline$kM\in\{N\}$&$eM\in\{M\}$\\
\hline$hL\in\{N\}$&$fL\in\{M\}$\\
\hline$Nk\in\{k\}$&$Mk\in\{e\}$\\
\hline$Lk\in\{g\}$&$ge\in\{g\}$\\
\hline$ke\in\{k\}$&$ee'\in\{e\}$\\
\hline$hg\in\{k\}$&$fg\in\{e\}$\\
\hline$Nh\in\{h\}$&$Mh\in\{f\}$\\
\hline$Lh=1$&$kf\in\{h\}$\\
\hline$gf=1$&$ef\in\{f\}$\\
\hline
\end{tabular}$$\\
\\
les \'equations d'associativit\'e sont donn\'ees par 8 lignes:\\
\begin{tabular}{|c|c|c|c|c|c|}
\hline $fgf\in\{f\}$&$fLh\in\{f\}$&$gfg\in\{g\}$&$fge\in\{e\}$&$fLN\in\{M\}$  \\
\hline $hgf\in\{h\}$&$hLh\in\{h\}$&$efg\in\{e\}$&$hge\in\{k\}$&$hLN\in\{N\}$\\
\hline $gef=1$&$gMh=1$&$kfg\in\{k\}$&$gee'\in\{g\}$&$gMN\in\{L\}$ \\
\hline $ee'f\in\{f\}$&$eMh\in\{f\}$&$Lhg\in\{g\}$&$ee'e"\in\{e\}$&$eMN\in\{M\}$ \\
\hline $kef\in\{h\}$&$kMh\in\{h\}$&$Mhg\in\{e\}$&$kee'\in\{k\}$&$kMN\in\{N\}$ \\
\hline $Lkf=1$&$LNh=1$&$Nhg\in\{k\}$&$Lke\in\{g\}$&$LNN'\in\{L\}$ \\
\hline $Mkf\in\{f\}$&$MNh\in\{f\}$&$Mke\in\{e\}$&$MNN'\in\{M\}$&$MNk\in\{e\}$ \\
\hline $Nkf\in\{h\}$&$NN'h\in\{h\}$&$Nke\in\{k\}$&$NN'N"\in\{N\}$&$NN'k\in\{k\}$ \\
\hline $MhL\in\{M\}$&$kfL\in\{N\}$&$NhL\in\{N\}$&$fgM\in\{M\}$&$geM\in\{L\}$\\
\hline $keM\in\{N\}$&$LkM\in\{L\}$&$NkM'\in\{M\}$&$NkM\in\{N\}$&$eMk\in\{e\}$\\
\hline$fkL\in\{e\}$&$hLk\in\{k\}$&$gMk\in\{g\}$ &$LhL'\in\{L\}$&$kMk'\in\{k\}$\\
\hline$gfL\in\{L\}$&$efL\in\{M\}$&$LNk\in\{g\}$&$ee'M\in\{M\}$&$hgM\in\{N\}$\\
\hline
\end{tabular}\\
\\
Nous revenons au but: il faut d\'emontrer que $m\geq bc$ et $q\geq ap$. \\
On suppose que $m<bc$ alors il y a 3 cas: \\
--il existe $h\neq h'$,$g\neq g'$ tels que $hg=h'g'$ alors
$L(hg)=L(h'g')$ donc $(Lh)g=(Lh')g'$ alors $g=g'$ car
($Lh=1$) et $L$ arbitraire; \\
--il existe $h$ ,$g\neq g'$ tels que $hg=hg'$ alors $L(hg)=L(hg')$ donc $(Lh)g=(Lh)g'$ alors $g=g'$ car ($Lh=1$); \\
--il existe $h\neq h'$, g tels que $hg=h'g$ alors $(hg)f=(h'g')f$ donc $h(gf)=h'(g'f)$ alors $h=h'$ car ($gf=1$)et $f$ arbitraire. \\
De m\^eme pour $q\geq ap$.\\
On pose $q< ap$ alors il existe 3 cas:\\
--il existe $f\neq f'$,$L\neq L'$ tels que $fL=f'L'$ alors
$g(fL)=g(f'L')$ donc $(gf)L=(gf')L'$ alors $L=L'$ car
($gf=1$)et $g$ arbitraire; \\
--il existe $f$,$L\neq L'$ tels que $fL=fL'$ alors $g(fL)=g(fL')$ donc $(gf)L=(gf)L'$ alors $L=L'$ car ($gf=1$) et $g$ arbitraire; \\
--il existe $f\neq f'$,L tels que $fL=f'L'$ alors $(fL)h=(f'L')h$ donc $f(Lh)=f'(Lh)$ alors $f=f'$ car ($Lh=1$) et $h$ arbitraire. \\
\\
Finalement :\\
\begin{theoreme}
Si $M$ est une matrice  strictement positive de taille $3$
telle que $m_{11}=1$ et $m_{22},m_{33}>1$, une condition
n\'ecessaire pour qu'elle marche est que $n\geq ac+1$, $r\geq bp+1$
et $m\geq bc$ , $q\geq ap$. 
\end{theoreme}

\noindent
\textbf{Rq:}
$ge=g$ et $e^2=e$ pour tout $e\in A(x_2,x_2)$ et $g\in A(x_2,x_1),f\in A(x_1,x_2)$ tel que$fg=e$ en effet:\\
$fg=e$ alors  $g(fg)=ge$ donc $(gf)g=g=ge$ aussi $fg=e$  alors $f(ge)=e$ donc $(fg)e=e$ alors $e^2=e$. \\
\\
\begin{theoreme}
Soit $M$ une  matrice triple d\'efinie par: \\
\begin{displaymath}
\mathbf{M} = \left( \begin{array}{ccc}
1 & a &b  \\
c & n &m\\
p & q &r
\end{array} \right)
\end{displaymath}
avec a,b,c,n,m,p,q,r$>1$.\\
alors $Cat(M)\neq{\Oe}\Leftarrow n=ac+1,r=bp+1,m=bc,q=ap$.\\
\end{theoreme}
\textbf{Preuve:} soit $A$ une semi-cat\'egorie dont les objets $\{x_1,x_2,x_3\}$ et les morphismes d\'efinis par:\\
$A(x_1,x_1)=\{1\}$,\\
$A(x_1,x_2)=\{f_1,...,f_a\}$,\\
$A(x_1,x_3)=\{h_1,...,h_b\}$,\\
$A(x_2,x_1)=\{g_1,...,g_c\}$,\\
$A(x_2,x_3)=\{k_1^1,...,k_c^b\}$,\\
$A(x_3,x_1)=\{L_1,...,L_p\}$,\\
$A(x_3,x_2)=\{M_1,...,M_p^a\}$,\\
$A(x_2,x_2)=\{e_1^1,...,e_c^a\}$ avec $e_i^j\neq 1$ pour tout i, \\
$A(x_3,x_3)=\{N_1^1,...,N_p^a\}$ avec $N_i^j\neq 1$ pour tout i. \\
\\
Les \'equations de la loi de composition sont d\'efinies par:
$$\begin{tabular}{|c|c|c|c|c|c|}
\hline$N_j^iN_{j'}^{i'}=N_{j'}^{i}$&$M_j^iN_{j'}^{i'}=M_{j'}^i$\\
\hline$L_{j'}N_j^i=L_j$&$g_{j'}M_j^i=L_j$\\
\hline$k_j^iM_{j'}^{i'}=N_{j'}^i$&$e_j^iM_{j'}^{i'}=M_{j'}^i$\\
\hline$h_iL_j=N_j^i$&$f_iL_j=M_j^i$\\
\hline$N_j^ik_{j'}^{i'}=k_{j'}^i$&$M_j^ik_{j'}^{i'}=e_{j'}^i$\\
\hline$L_{i'}k_j^i=g_j$&$g_{j'}e_j^i=g_j$\\
\hline$k_{j'}^{i'}e_j^i=k_{j}^{i'}$&$e_j^ie_{j'}^{i'}=e_{j'}^{i}$\\
\hline$h_ig_j=k_j^i$&$f_ig_j=e_j^i$\\
\hline$N_j^ih_{j'}=h_i$&$M_j^ih_{j'}=f_i$\\
\hline$L_ih_j=1$&$K^i_jf_{j'}=h_i$\\
\hline$g_if_j=1$&$e_j^if_{j'}=f_i$\\
\hline
\end{tabular}$$\\
\\
Toutes les \'equations associatives marchent. Donc $A$ est une semi-cat\'egorie. \\
Soit $B=A\oplus \{id_{x_2}\}\oplus\{id_{x_3}\}$ donc $B$ est une cat\'egorie associ\'ee \`a $M$ ce qui donne $Cat(M)\neq {\Oe} $. \\
\\
\textbf{Notation:}\\
On a pour la matrice triple $M$ d\'efinie par:
\begin{displaymath}
\mathbf{M} = \left( \begin{array}{ccc}
1 & a &b  \\
c & n &m\\
p & q &r
\end{array} \right)
\end{displaymath}
avec $n=ac$ ,$r=bp$, $m=bc$ et $q=ap$ d'apres ce qui pr\'ec\`ede $M$ admet $A$ comme semi-cat\'egorie. \\
\\
Maintenant on va chercher une semi-cat\'egorie associ\'ee \`a $M$ avec $n>ac,r>bp,m>bc , q>ap $ et apr\`es on ajoutera les identit\'es. \\
Soit $M(ac+1)$ matrice d\'efinie par :\\
\begin{displaymath}
\mathbf{M(ac+1)} = \left( \begin{array}{ccc}
1 & a    &b  \\
c & ac+1 &m\\
p & q    &r
\end{array} \right)
\end{displaymath}
avec $r=bp,m=bc ,q=ap $.\\
Soit $A'$ une semi-cat'egorie dont les objets sont $Ob(A')=Ob(A)$  avec les morphismes: \\
$A'(x_1,x_1)=\{1\}$;\\
$A'(x_1,x_2)=\{f_1,...,f_a\}$;\\
$A'(x_1,x_3)=\{h_1,...,h_b\}$;\\
$A'(x_2,x_1)=\{g_1,...,g_c\}$;\\
$A'(x_2,x_3)=\{k_1^1,...,k_c^b\}$;\\
$A'(x_3,x_1)=\{L_1,...,L_p\}$;\\
$A'(x_3,x_2)=\{M_1,...,M_p^a\}$;\\
$A'(x_2,x_2)=\{e_1^1,...,e_c^a\}\cup\{e'\}$ avec $e'\not\in A(x_2,x_2)$; \\
$A'(x_3,x_3)=\{N_1^1,...,N_p^a\}$  .\\
Les \'equations de la loi de composition de $A'$ sont les m\^emes que pour $A$, en plus les \'equations dependant de $e'$ sont:\\
$e'f_i=e^1_1f_i$;\\
$e'M_{j}^{i}=e^1_1 M_{j}^{i}$;\\
$e'e^i_j=e^1_1 e^i_j$;\\
$e^i_j e'=e^i_je^a_c $;\\
$k^i_j e'=k^i_j e^a_c$;\\
$g_i e'=g_i e^a_c$.\\
Les \'equations associatives marchent mais dans le cas o\`u il y a $e'$ les \'equations sont:\\
$$\begin{tabular}{|c|c|c|c|c|}
\hline $g_ne'M^i_j$&$e'M^i_jk^l_n$&$e'f_nL_i$&$e^i_je'M^l_n$\\
\hline$ee'e"$&$k^i_jee'$&$L_nk^i_je'$&$e'M^i_jN^l_n$\\
\hline$N_j^ik^l_ne'$&$f_ig_je'$&$h_ig_je'$&$g_le^i_je'$\\
\hline$k^i_je'M^l_n$&$e'M^i_jh_n$&$e'f_ig_j$&$M_j^ik^l_ne'$\\
\hline$g_je'f_i$&$e^i_j e' f_{i'}$&$k^i_je'f_{i'}$&$k^i_je'f_{i'}$\\
\hline
\end{tabular}$$\\
on va v\'erifier ces \'equations. \\
$g_j(e'f_i)=g_j(e^1_1 f_i)=g_jf_1=1=(g_je')f_i=(g_je^a_c)f_i=gf=1$ vraie.\\
$(e^i_j e') f_{i'}=(e^i_j e^1_1)f_{i'}=e^i_1f_{i'}=f_i=e^i_j(e' f_{i'})=e^i_j f_1=f_i$ vraie de meme $e'ef$.\\
$(k^i_je')f_{i'}=(k^i_c)f_{i'}=h_i=k^i_j(e'f_{i'})=k^i_jf_1=h_i$ varie .\\
$(k^i_je')M^l_n=k^i_cM^l_n=N^i_1=k^i_j(e'M^l_n)=k^i_jM^1_n=N^i_1$ varie.\\
$(e'M^i_j)h_n=M^1_jh_n=f_1=e'(M^i_jh_n)=e'f_i=f_1$ vraie.\\
$(e'f_i)g_j=f_1g_j=e^1_j=e'(f_ig_j)=e'e^i_j=e^1_j$ vraie;\\
$(M_j^ik^l_n)e'=e^i_ne'=e^i_c=M_j^i(k^l_ne')=M_j^iK^1_c=e^i_c$ vraie.\\
$(N_j^ik^l_n)e'=k^i_ne'=k^i_c=N_j^i(k^l_ne')=N_j^ik^l_c=k^i_c$ vraie.\\
$(f_ig_j)e'=e^i_je'=e^i_c=f_i(g_je')=f_ig_c=e^i_c$ vraie.\\
$(h_ig_j)e'=k^i_je'=k^i_c=h_i(g_je')=h_ig_c=k^i_c$ vraie.\\
$(g_le^i_j)e'=g_je'=g_c=g_l(e^i_je')=g_le^i_c=g_c$ vraie de meme $ge'e$.\\
$(e^i_je')e^l_n=e^i_ce^l_n=e^i_n=e^i_j(e'e^l_n)=e^i_je^1_n=e^i_n$ vraie de meme $e'ee"$ et $ee"e'$.\\
$(k^l_ne^i_j)e'=k^l_je'=k^l_c=k^l_n(e^i_je')=k^l_ne^i_c=k^l_c$ vraie de meme $ke'e$.\\
$(L_nk^i_j)e'=g_je'=g_c=L_n(k^i_je')=L_nk^i_c=g_c$ vraie .\\
$(e'M^i_j)N^l_n=M^1_jN^l_n=M^1_n=e'(M^i_jN^l_n)=e'M^i_n=M^1_n$ vraie.\\
$(g_ne')M^i_j=g_cM^i_j=Lj=g_n(e'M^i_j)=g_nM^1_j=L_j$ vraie .\\
$(e'M^i_j)k^l_n=M^1_jk^l_n=e^1_n=e'(M^i_jk^l_n)=e'e^i_n=e^1_n$ vraie.\\
$(e'f_n)L_i=f_1L_i=M^1_i=e'(f_nL_i)=e'M^n_i=M^1_i$ vraie.\\
$(e^i_je')M^l_n=e^i_cM^l_n=M^i_n=e^i_j(e'M^l_n)=M^i_n$ vraie de meme $e'eM$.\\
Donc $A_{1}$ est une semi-cat\'egorie, on note $e'=e_1$ \\
on suppose que $A_{n-1}$ est une semi-cat\'egorie telle que $A_{n-1}=A'\cup\{e_2,...,e_{n-1}\}=A\cup \{e_1,...,e_{n-1}\}$ avec $e_i\neq e_j$ pour tout $i,j\in \{1,...,(n-1)\}$  \\
avec la loi de composition d\'efinie  par:\\
$e_i(...)=e_1^1(...)$,\\
$(...)e_i=(...)e_c^a$,\\
$e_ie_j=e^1_c=e_je_i$,\\
${e_i}^2=e_i$.\\
Pour v\'erifier que $A_{n-1}$ est une semi-cat\'egorie on consid\`ere les \'equations de la loi d'associativit\'e suivantes.\\
$(e_ie_j)f_v=e_c^1f_v=f_1=e_i(e_jf_v)=e_if_1=f_1$ vraie.\\
$g_v(e_ie_j)=g_ve_c^1=g_c=(g_ve_i)e_j=g_ce_j=g_c$ vraie .\\
$(e_ie_j)M_s^t=e_c^1M_s^t=M_s^1=e_i(e_jM_s^t)=e_iM_s^1=M_s^1$ vraie.\\
$(k_s^te_i)e_j=k_c^te_j=k^t_c=k_s^t(e_ie_j)=k_s^te_c^1=k^t_c$ vraie.\\
$(e_ie_j)e_k=e_c^1e_k=e_c^1e_c^a=e^1_c=e_i(e_je_k)=e_i(e_c^1)=e_1^1e_c^1=e^1_c$ varie. \\
$(e_j^ie_k)e_l=(e_j^ie_c^a)e_l=e_c^ie_c^a=e^i_c=e_j^ie_c^1=e_c^i$ vraie .\\
$(e_ke_j^i)e_l=(e_1^1e_j^i)e_l=e_j^1e_c^a=e^1_c=e_1^1(e_j^ie_c^a)=e_c^1$ vraie. \\
Les autres \'equations ressemblent aux \'equations dans le cas de $e_1$. \\
Donc $A_{n-1}$ est une semi-cat\'egorie associe\'e \`a la matrice suivant:\\
\begin{displaymath}
\mathbf{M(ac+(n-1))} = \left( \begin{array}{ccc}
1 & a        &b  \\
c & ac+(n-1) &bc\\
p & ap       &bp
\end{array} \right)
\end{displaymath}
Maintenant on ajoute des morphismes sur  $A(x_2,x_3)$. \\
Soient $k_1,...,k_m$ morphismes dans $A(x_2,x_3)$  avec $k_i\neq k_j$ pour tout $i,j\in \{1,...,m\}$ tel que le loi de composition est d\'efinie par:\\
$k_i(...)=k_1^1$\\
$(...)k_i=(...)k_c^b$\\
$k_je_i=k_1^1e_c^a=k_c^1$.\\
Les \'equations associatives associ\'ees \`a $k_i$ sont:\\
$(k_ve_j^i)f_l=(k^1_1e_j^i)f_l=k_j^1f_l=h_1=k_v(e_j^if_l)=k_vf_i=k_1^1f_i=h_1$ vraie.\\
$(k_ve_i)f_l=(k^1_1e_c^a)f_l=k_c^1f_l=h_1=k_v(e_if_l)=k_vf_1=k_1^1f_1=h_1$ vraie.\\
$(k_vM_j^i)h_l=(k^1_1M_j^i)h_l=N_j^1h_l=h_1=k_v(M_j^ih_l)=k_vf_i=k_1^1f_i=h_1$ vraie.\\
$(k_vM_j^i)N_l^n=(k^1_1M_j^i)N_l^n=N_j^1N_l^n=N_l^1=k_v(M_j^iN_l^n)=k_vM_l^i=N_l^1$ vraie.\\
$(k_ne_j^i)M_v^o=(k^1_1e_j^i)M_v^o=k_j^1M_v^o=N_v^1=k_n(e_j^iM_v^o)=k_nM_v^i=N_v^1$ vraie.\\
$(k_ne_j)M_v^o=(k^1_1e_c^a)M_v^o=k_c^1M_v^o=N_v^1=k_n(e_jM_v^o)=k_nM_v^1=N_v^1$ vraie.\\
$(k_vM_j^i)k_d^o=(k^1_1M_j^i)k_d^o=N_j^1k_d^o=k_d^1=k_v(M_j^ik_d^o)=k_ve_d^i=k_d^1$ vraie.\\
$(k_vM_j^i)k_d=(k^1_1M_j^i)k_d=N_j^1k_c^b=k_c^1=k_v(M_j^ik_d)=k_ve_c^i=k_c^1$ vraie.\\
$(k_d^oM_j^i)k_v=(N_j^o)k_v=N_j^ok_c^b=k_c^o=k_d^o(M_j^ik_v)=k_d^oe_c^i=k_c^o$ vraie.\\
$(k_vf_i)L_j=(k_1^1f_i)L_j=h_1L_j=N_j^1=k_v(f_iL_j)=k_1^1M_j^i=N_j^1$ vraie.\\
$(k_vf_i)g_j=(k_1^1f_i)g_j=h_1g_j=k_j^1=k_v(f_ig_j)=k_1^1e_j^i=k_j^1$ vraie.\\
$(L_jk_v)f_i=(L_jk_c^b)f_i=g_cf_i=1=L_j(k_vf_i)=1$ vraie.\\
$(L_jk_i)M_v^o=(L_jk_c^b)M_v^o=g_cM_v^o=L_v=L_j(k_iM_v^o)=L_jN_v^1=L_v$ vraie.\\
$(M_j^ik_o)f_v=(M_j^ik_c^b)f_v=e_c^if_v=f_i=M_j^i(k_if_v)=M_j^ih_1=f_i$ vraie.\\
$(L_jk_i)e_v^o=(L_jk_c^b)e_v^o=g_ce_v^o=g_v=L_j(k_ie_v^o)=L_jk_v^1=g_v$ vraie.\\
$(N_j^ik_o)f_v=(N_j^ik_c^b)f_v=k_c^if_v=h_i=N_j^i(k_if_v)=N_j^ih_1=h_i$ vraie.\\
$(M_j^ik_o)e_v^n=(M_j^ik_c^b)e_v^n=e_c^ie_v^n=e_v^i=M_j^i(k_ie_v^n)=M_j^ik_v^1=e_v^i$ vraie.\\
$(M_j^ik_o)e_v=(M_j^ik_c^b)e_v=e_c^ie_v=e_c^i=M_j^i(k_oe_v)=M_j^ik_c^1=e_c^i$ vraie.\\
$(M_j^ik_o)M_v^n=(M_j^ik_c^b)M_v^n=e_c^iM_v^n=M_v^i=M_j^i(k_iM_v^n)=M_j^iN_v^1=M_v^i$ vraie.\\
$(N_j^ik_o)e_v^n=(N_j^ik_c^b)e_v^n=k_c^ie_v^n=k_v^i=N_j^i(k_ie_v^n)=N_j^ik_v^1=k_v^i$ vraie.\\
$(N_j^ik_o)e_v=(N_j^ik_c^b)e_v=k_c^ie_c^a=k_c^i=N_j^i(k_ie_v^n)=N_j^ik_c^1=k_c^i$ vraie.\\
$(N_j^ik_o)M_v^n=(N_j^ik_c^b)M_v^n=k_c^iM_v^n=N_v^i=N_j^i(k_oM_v^n)=N_j^iN_v^1=N_v^i$ vraie.\\
$(f_iL_j)k_o=M_j^ik_o=e_c^i=f_i(L_jk_o)=f_ig_c=e_c^i$ vraie.\\
$(h_iL_j)k_o=N_j^ik_o=k_c^i=h_i(L_jk_o)=h_ig_c=k_c^i$ vraie.\\
$(g_vM_j^i)k_o=L_jk_o=g_c=g_v(M_j^ik_o)=g_ve_c^i=g_c$ vraie.\\
$(e_v^nM_j^i)k_o=M_j^nk_o=e_c^n=e_v^n(M_j^ik_o)=e_v^ne_c^i=e_c^n$ vraie.\\
$(e_vM_j^i)k_o=M_j^1k_o=e_c^1=e_v(M_j^ik_o)=e_ve_c^i=e_c^1$ vraie.\\
$(L_vN_j^i)k_o=L_jk_o=g_c=L_v(N_j^ik_o)=L_vk_c^i=g_c$ vraie .\\
$(M_v^nN_j^i)k_o=M_j^nk_o=e_c^n=M_v^n(N_j^ik_o)=M_v^nk_c^i=e_c^n$ vraie .\\
$(N_v^nN_j^i)k_o=N_j^nk_o=k_c^n=N_v^n(N_j^ik_o)=N_v^nk_c^i=k_c^n$ vraie .\\
\\
Par la meme construction on peux ajouter aussi sur $A(x_3,x_3)$ et sur $A(x_3,x_2)$ des morphismes adjoints avec la d\'efinition de la loi de composition:\\
$N_i(...)=N_1^1(...)$ .\\
$(...)N_i=(...)N_p^b$.\\
$N_iN_j=N^1_p=N_jN_i$.\\
$N_iK_j=K_c^1$ et $M_iN_j=M_p^1$.\\
${N_i}^2=N_i$ pour tout $i\in \{1,...,r-1\}$.\\
$M_i(...)=M_1^1(...)$ .\\
$(...)M_i=(...)M_P^a$.\\
$M_ik_j=e_c^1$ ,$e_iM_j=M_p^1$ et $k_iM_j=N_p^1$. \\
pour tout $i\in \{1,...,q\}$.\\
Les \'equations associatives marchent, donc $A_{(n-1,m,q,r-1)}$ est une semi-cat\'egorie  \\
c.\`a .d  en ajoutant les identites $A'_{(n,m,q,r)}$ est une cat\'egorie  associ\'ee \`a la matrice $M$ qui est d\'efinie par:\\
\begin{displaymath}
\mathbf{M} = \left( \begin{array}{ccc}
1 & a        &b\\
c & ac+n     &bc+m\\
p & ap+q       &bp+r
\end{array} \right)
\end{displaymath}
o\`u $n,r,m,q$ sont des entiers naturels.\\
\begin{corollary}
Soit $M$ une matrice d'order 3 tel que :
\begin{displaymath}
\mathbf{M} = \left( \begin{array}{ccc}
z & a        &b\\
c & n     &m\\
p & q       &r
\end{array} \right)
\end{displaymath}\\
avec $z\geq 1$,$n>ac$ ,$r>bp$, $m\geq bc$ et $q \geq ap$ alors $Cat(M)\neq {\Oe}$\\
\end{corollary}
\textbf{Preuve}: soit $N$ une matrice d\'efinie par:
\begin{displaymath}
\mathbf{N} = \left( \begin{array}{ccc}
1 & a        &b\\
c & n     &m\\
p & q       &r
\end{array} \right) .
\end{displaymath}\\
D'apr\'es le th\'eor\`eme pr\'ec\'edant $Cat(N)\neq{\Oe}$ alors il
existe une cat\'egorie $A$ d\'efinie comme pr\'ec\'edemment, soit $A'$
une cat\'egorie dont les objets $Ob(A')=Ob(A)$ et les morphismes
$Mor(A')=Mor(A)\cup\{n_1,n_2,...,n_{z-1}\}$ avec
$A'(x_1,x_1)=\{1,n_1,n_2,...,n_{z-1}\}$ et les \'equations de la loi
de composition associ\'ees \`a $n_i$ d\'efinies par:\\
$n_i(...)=(...)$\\
$(...)n_i=(...)$\\
$n_i^2=n_i$\\
$n_in_j=n_1$ avec $i\neq j$. \\
Alors $A'$ une cat\'egorie associ\'ee \`a $M$ donc
$Cat(M)\neq{\Oe}$.\\
\\
\section {Matrices G\'en\'erales}
\begin{theoreme}
\label{gen1}
Soit $M$ une matrice de taille $n$ telle que $M$ d\'efinie par:\\
\begin{displaymath}
\mathbf{M} = \left( \begin{array}{cccc}
1 & M_{12} & \ldots &M_{1n}\\
M_{21} & M_{22} & \ldots & M_{2n}\\
\vdots & \vdots & \ddots&\vdots\\
M_{n1} & M_{n2} &  \ldots   &M_{nn}
\end{array} \right)
\end{displaymath}
avec $M_{ij}> 0$ $ \forall i,j \in\{1,...,n\}$ et $M_{ii}>1 $ pour $\in\{2,...,n\}$,\\
alors $Cat(M)\neq {\Oe}$ si et seulement si $M_{ii}> M_{1i}M_{i1}
\forall i \in\{2,...,n\}$ et $M_{ij}\geq M_{i1}M_{1j}$ avec $i,j
\in\{2,...,n\}$.\\
\end{theoreme}
En effet:\\
on supose que $M$ marche alors il existe $A$ cat\'egorie
associ\'ee \`a $M$ dont les objets sont $\{x_1,...,x_n\}$ et
$|A(x_i,x_j)|=M_{ij}$.\\
On va d\'emontrer que $M_{ii}> M_{1i}M_{i1}$ on suppose que
$M_{ii}\leq M_{1i}M_{i1}$, soit $A(x_1,x_i)=\{f_1,...,f_a\}$ et
$A(x_i,x_1)=\{g_1,...,g_b\}$ et $A(x_i,x_i)=\{1,e_2,...,e_c\}$
avec $a=M_{1i}$ , $b=M_{i1}$ et $c=M_{ii}$ on a $gf=1$ pour tout
f,g et $fg=e$ on a $M_{ii}\leq M_{1i}M_{i1}$  alors soit il existe
$g,g',f,e$ tel que $fg=fg'=e$ alors $g(fg)=g(fg')=ge$ donc
$(gf)g=(gf)g'$ alors $g=g'$ impossible car $g\neq g'$; soit de
m\^eme pour si on a $f,f'$ tel que $fg=f'g=e$ impossible, soit si
on a $f,f',g,g'$ avec $fg=f'g'=e$ impossible aussi, ce qui donne
$M_{ii}\geq M_{1i}M_{i1}$. Si $\exists f,g $ tel que $fg=1$ alors
$g'(fg)=g'$
alors $(g'f)g=g'$ donc $g=g'$ impossible.\\
Finalement $M_{ii}> M_{1i}M_{i1}$.\\
Pour $M_{ij}\geq M_{i1}M_{1j}$,\\
soit $A(x_i,x_1)=\{g_1,...,g_b\}$,$A(x_1,x_j)=\{h_1,...,h_m\}$ et
$A(x_i,x_j)=\{L_1,...,L_v\}$ avec $b=M_{i1},m=M_{1j}$ et
$v=M_{ij}$, on suppose que $M_{ij}< M_{i1}M_{1j}$ , alors ils
$\exists L,h,h',g $ou $L,h,h',g,g'$ ou $L,h,g,g'$ les trois sont
les m\^emes type de d\'emonstration; je veux prendre le cas o\`u
$\exists L,h,h',g,g'$ tel que $L=hg=h'g'$ alors $Lf=(hg)f=(h'g')f$
donc $Lf=h(gf)=h'(g'f)$ alors $Lf=h=h'$ impossible car $h\neq h'$,
ce qui donne $M_{ij}\geq M_{i1}M_{1j}$ pour tout $i,j
\in\{2,...,n\}$.\\
\\
Maintenant on va d\'emontrer le sens inverse: si $M_{ii}=
M_{1i}M_{i1} \forall i \in\{2,...,n\}$ et $M_{ij}=
M_{i1}M_{1j}$ avec $i,j \in\{2,...,n\}$ alors $M$ marche .\\
En effet:\\
Soit $A'$ une semi-cat\'egorie dont les objets sont
$\{x_1,...,x_n\}$ et
$|A'(x_i,x_j)|=M_{ij}$ avec les notations suivantes:\\
\\
$A'(x_1,x_1)= 1$.\\
\\
$A'(x_1,x_i)=\{{_if_1},...,{_if_{M_{1i}}}\}$ pour tout $i\in\{2,...,n\}$.\\
\\
$A'(x_i,x_1)=\{{_ig_1},...,{_ig_{M_{i1}}}\}$ pour tout $i\in\{2,...,n\}$.\\
\\
$A'(x_i,x_i)=\{{_ie_1^1},...,{_ie^{M_{1i}}_{M_{i1}}}\}$  pour tout $i\in\{2,...,n\}$.\\
\\
$A'(x_i,x_j)=\{{_j^iH_1^1},...,{_j^iH^{M_{1j}}_{M_{i1}}}\}$ pour tout $i,j\in\{2,...,n\}$ avec $i \neq j$. \\
Les \'equations de la loi de composition sont d\'efinies par:\\
$$\begin{tabular}{|c|c|c|c|c|c|}
\hline$_i^jH_b^a \circ{_j^iH_d^c}={_ie^a_d}$&$_je^a_b\circ{_j^iH_d^c}={_j^iH_d^a}$&$_if_a \circ{ _jg_b}={_i^jH_b^a}$\\
\hline$_ig_c \circ {_ie_b^a}={_ig_b}$& $_ie_b^a\circ{_ie^c_d}={_ie^a_d}$& $_ie_b^a \circ{ _ie_b^a}={_ie_b^a}$\\
\hline $g\circ f=1$& $_ie_b^a \circ {_if_c}={_if_a}$& $_j^iH_b^a\circ{_jf_c}={_if_a}$\\
\hline$_ig_a\circ{ _i^jH_c^b}={_jg_c}$&$_i^jH^a_b \circ{ _je_d^c}={_i^jH_d^a}$&$_i^jH_b^a \circ  {_{j'}^{i'}H_d^c}={_i^{i'}H^a_d}$\\
\hline$_ig_a\circ {_j^iH_c^b}={_ig_c}$& $_if_a \circ {_ig_b}={_ie^a_b}$&$_jf_a \circ{ _ig_b}={_j^iH^a_b}$\\
\hline
\end{tabular}$$\\
\\
D'apr\`es cette d\'efinition les \'equations associatives
marchent comme dans l'exemple de la matrice triple.\\
Alors $A'$ est une semi-cat\'egorie associ\'ee \`a $M'$ tel que :
\begin{displaymath}
\mathbf{M'} = \left( \begin{array}{cccc}
1 & M_{12} & \ldots &M_{1n}\\
M_{21} & (M_{21} M_{12}) & \ldots & (M_{21} M_{1n})\\
\vdots & \vdots & \ddots&\vdots\\
M_{n1} & (M_{n1}M_{12}) &  \ldots   &(M_{n1}M_{1n})
\end{array} \right) .
\end{displaymath}
On ajoute des morphismes pour g\'en\`eraliser le th\'eor\`eme sur
les matrices de taille $3$. On arrive surtout aux ensembles des morphismes suivants:\\
\\
$A''(x_i,x_i)=\{_ie_1^1,...,{_i{e^{M_{1i}}_{M_{i1}}},{_ie_1},...,{_ie_{s_i}}\}}$ pour tout $i\in\{2,...,n\}$\\
\\
$A''(x_i,x_j)=\{{_j^iH_1^1},...,{_j^i{H^{M_{1j}}_{M_{i1}}},{_j^iH_1},...,{_j^iH_{t_j^i}}\}}$ pour tout $i,j\in\{2,...,n\}$ avec $i \neq j$ \\
\`a condition que tous les ajout\'es sont distincts, et en plus
la loi de composition est d\'efinie par :\\
\\
${_ie_k}\circ(...)={_ie_1}\circ (...)$ pour tout $k\in\{1,...,s_i\}$ et $i\in\{2,...,n\}$. \\
\\
$(...)\circ{_ie_k}=(...) \circ {_ie_{M_{i1}}^{M_{1i}}} $ pour tout $k\in\{1,...,s_i\}$ et $i\in\{2,...,n\}$. \\
\\
${_ie_k}\circ {_ie_k}={_ie_k}$ pour tout $k\in\{1,...,s_i\}$ et $i\in\{2,...,n\}$. \\
\\
${_ie_k}\circ {_ie_p}={_ie^1_{M_{i1}}}$ pour tout $k,p\in\{1,...,s_i\}$ et $i\in\{2,...,n\}$.\\
\\
${_j^iH_p}\circ (...)={_j^iH_1^1}\circ (...)$ pour tout $p\in\{1,...,t_j^i\}$ et $i,j\in\{2,...,n\}$.\\
\\
$(...)\circ {_j^iH_p}=(...)\circ _jH_{M_{i1}}^{M{1j}}$ pour tout  $p\in\{1,...,t_j^i\}$ et $i,j\in\{2,...,n\}$.\\
\\
${_j^iH_p}\circ {_ie_k}={_j^iH^1_{M_{i1}}}$ pour tout $k\in\{1,...,s_i\}$,$p\in\{1,...,t_j^i\}$ et $i,j\in\{2,...,n\}$.\\
\\
${_je_k}\circ{_j^iH_p}={_j^iH^1_{M_{i1}}}$  pour tout $k\in\{1,...,s_i\}$,$p\in\{1,...,t_j^i\}$ et $i,j\in\{2,...,n\}$.\\
\\
Donc $A''$ une semi-cat\'gorie associ\'ee \`a la matrice $M''$ d\'efinie par:\\
\begin{displaymath}
\mathbf{M''} = \left( \begin{array}{cccc}
1 & M_{12} & \ldots &M_{1n}\\
M_{21} & (M_{21} M_{12})+s_i & \ldots & (M_{21} M_{1n})+t_n^2\\
\vdots & \vdots & \ddots&\vdots\\
M_{n1} & (M_{n1}M_{12}+)+t_2^n &  \ldots   &(M_{n1}M_{1n})+s_n
\end{array} \right) .
\end{displaymath}
On peut ensuite ajouter les identit\'es sur $x_2,\ldots , x_n$. \\
\textbf{Finalement} si $M=(M_{ij})_n$ une matrice positive d'order
$n$ telle que $M_{11}=1$ alors $Cat(M)\neq{\Oe}$ si et seulement
si $M_{ii}> M_{1i}M_{i1} \forall i \in\{1,...,n\}$ et $M_{ij}\geq
M_{i1}M_{1j} \forall i,j \in\{1,...,n\} i\neq j$.\\
On obtient le th\'eor\`eme suivant:

\begin{theoreme}
Si $M_{11}=1$ et $M_{ii}>1$ pour $i>1$, avec $M_{ij}>0$ $\forall i,j$, alors 
$Cat(M)\neq \Oe$  si et seulement
 si $M_{ii}> M_{1i}M_{i1}$ $\forall i \in\{1,...,n\}$ et $M_{ij}\geq
M_{i1}M_{1j}$ $\forall i\neq j \in\{1,...,n\}$.\\
\end{theoreme}

On doit maintenant traiter la possibilit\'e que $M_{ii}=1$ pour plusieurs  $i$ distincts. 

\begin{definition}:
Soit $A$ une categorie d'ordre $n$ avec objets $x_1,\ldots , x_n$,
on dit que $x_i$ et $x_j$ sont {\em isomorphes} s'il existe $f\in
A(x_i,x_j)$ et $g\in A(x_j,x_i)$ tels que $fg = 1_{x_j}$ et $gf =
1_{x_i}$.
\end{definition}
 \textbf{Rq}: Si $x_i$ et $x_j$ sont isomorphes, alors pour tout objet $x_k$ on a des isomorphismes
d'ensembles
$$
A(x_k,x_i)\stackrel{\cong}{\rightarrow} A(x_k,x_j),
$$
donn\'es par $h\mapsto fh$ dans une direction, et $u\mapsto gu$
dans l'autre; et
$$
A(x_i,x_k)\stackrel{\cong}{\rightarrow} A(x_j,x_k),
$$
donn\'e par $h\mapsto hg$ dans une direction, et $u\mapsto uf$
dans l'autre. Si $M$ est la matrice de $A$, on en d\'eduit:
$$
\forall k, \;\; M_{ki} = M_{kj}$$et$$\forall k, \;\; M_{ik} =
M_{jk}.
$$ \\
\begin{definition}:
\label{defreduite} 
Soit $A$  une cat\'egorie telle qu'il existe deux
objets distincts $x_i$ et $x_j$ ($i\neq j$) qui sont isomorphes,
on dira que $A$ est {\em non-r\'eduite}. On dira que $A$ est {\em
r\'eduite} sinon,c'est-\`a-dire si deux objets distincts sont
toujours non-isomorphes. On dira qu'une matrice $M$ est {\em
non-r\'eduite} s'il existe $i\neq j$ tel que
$$\forall k, \;\;
M_{ki} = M_{kj}$$
et
$$\forall k, \;\; M_{ik} = M_{jk},$$
cela veut
dire que la ligne $i$ \'egale la ligne $j$ et la colonne $i$
\'egale la colonne $j$. On dira qu'une matrice $M$ est {\em
r\'eduite} si elle n'est pas non-r\'eduite.
\end{definition}
\textbf{Rq:} D'apres le debut ci-dessus, on obtient que si $A$ est
non-r\'eduite, alors $M$ est non-r\'eduite. Donc, par contrapos\'e
si $M$ est r\'eduite alors $A$ est r\'eduite. Le contraire n'est
pas forc\'emment vrai: il peut exister une cat\'egorie $A$ telle
que $M$ est non-r\'eduite, mais $A$ r\'eduite, par exemple on peut
avoir une cat\'egorie $A$ d'ordre $2$ dont la matrice
non-r\'eduite est
$$
M=\left(
\begin{array}{cc}2 & 2 \\2 & 2\end{array}\right)
$$
mais telle que les deux objets de $A$ sont non-isomorphes et donc
$A$ r\'eduite.
\begin{theoreme}
\label{reduire}
Si $M$ une matrice non r\'eduite, on peut r\'eduire
$M$ en une sous matrice $N$ r\'eduite telle que M marche si et
seulement si  $N$ marche.
\end{theoreme}
En effet: Supposons que $M$ est une matrice $n\times n$
non-r\'eduite. On peut d\'efinir une r\'elation d'\'equivalence
sur l'ensemble d'indices $\{ 1,\ldots ,n\}$ en disant que $i\sim
j$ si $\forall k, \;\; M_{ki} = M_{kj}$ et$\forall k, \;\; M_{ik}
= M_{jk}$. Celle-ci est sym\'etrique, reflexive et transitive. On
obtient donc une partitionde l'ensemble d'indices en r\'eunion
disjointe de sous-ensembles$$\{ 1,\ldots ,n\} = U_1\sqcup U_2
\sqcup \cdots \sqcup U_m$$avec $U_a\cap U_b= \emptyset$, telle que
tous les \'el\'ements d'un $U_a$ donn\'e sont \'equivalents, et
les \'el\'ements de $U_a$ ne sont pas \'equivalents aux
\'el\'ements de $U_b$ pour $a\neq b$. (Ce sont les classes
d'\'equivalence pour la r\'elation d'\'equivalence). Choisissons
un repr\'esantant $r(a)\in U_a$ pour chaque classe
d'\'equivalence. Dans l'autre sens, on note par $c(i)\in \{
1,\ldots , m\}$ l'unique \'el\'ement telle que $i\in U_{c(i)}$.
Ici $c(i)$ est la classe d'\'equivalence contenant $i$. On
a$$c(r(a))=a$$mais $r(c(i))$ n'est pas toujours \'egale \`a $i$:
on a seulement qu'ils sont \'equivalents $r(c(i))\sim i$.  On
obtient une sous-matrice de taille $m\times m$
$$
N_{ab}:= M_{r(a),r(b)}.
$$
On peut faire en sorte que $r(a)<r(b)$ pour
$a<b$: on choisit $r(a)$ le plus petit \'el\'ement de $U_a$, et on
num\'erote les classes $U_a$ par ordre croissant de leur plus
petit \'el\'ement. Dans ce cas $N$ est vraiement une sous-matrice
de $M$. On a $N$ r\'eduite, puisque les \'el\'ements de
$U_a$ et $U_b$ ne sont pas \'equivalents pour $a\neq b$. Si $A$
est une cat\'egorie dont la matrice est $M$, on obtient une
sous-cat\'egorie pleine $B\subset A$ qui consiste des objets
$r(a)$ seulement, $a=1,\ldots , m$. La matrice de $B$ est $N$. 
Donc si $M$ marche, alors $N$ marche. L'\'equivalence entre
$i$ et $r(c(i))$ implique que pour tout $k$ on a$$M_{k,i} =
M_{k,r(c(i))},\;\;\; M_{i,k} = M_{r(c(i)),k}.$$ On en d\'eduit que
pour tout $i,j$ on a $$M_{i,j}= M_{r(c(i)),j} = M_{r(c(i)),r(c(j))}
= N_{c(i),c(j)}.$$ Ceci indique comment aller dans l'autre
sens. Supposons que $B$ est une cat\'egorie dont la matrice est
$N$. Notons par $y_1,\ldots , y_m$ les objets de $B$. On d\'efinit
une cat\'egorie $A$ avec objets not\'es $x_1,\ldots , x_n$ en
posant$$A(x_i,x_j) \cong B(y_{c(i)},y_{c(j)}).$$ On pourrait d\'efinir $$A(x_i,x_j) := \{ (i,j,\beta),
\;\;\; \beta \in B(y_{c(i)},y_{c(j)}) \} .$$ La composition est la
m\^eme que celle de $B$, i.e.$$(i,j,\beta )(j,k,\beta '):=
(i,k,\beta \beta ').$$ De m\^eme pour les identit\'es, et les \'equations associatives
et les r\`egles des identit\'es sont faciles \`a
v\'erifier. Donc $A$ est une cat\'egorie.\\
On a:
$$
|A(x_i,x_j)| = |B(y_{c(i)},y_{c(j)})| = N_{c(i),c(j)} = M_{i,j}.
$$
Donc $A$ corr\'espond \`a la matrice $M$.\\ 
Finalement: \'etant donn\'ee une matrice non-r\'eduite $M$, on peut
construire par la construction pr\'ec\'edante une sous-matrice $N$
qui est r\'eduite, telle que $M$ marche si et seulement si $N$
marche. La sous-matrice $N$ est unique \`a permutation d'indices
pr\`es.

\begin{lemma}
\label{marchepas}
Soit $M$ est une matrice r\'eduite avec $M_{i,j}>0$, et s'il existe
$i\neq j$ tels que $M_{i,i}=1$ et $M_{j,j}=1$, alors $M$ ne marche
pas.
\end{lemma}
En effet: On suppose que $M$ marche alors il existe une
cat\'egorie $A$ associ\'ee \`a $M$ et comme $M$ est r\'eduite
alors $A$ est r\'eduite. En plus $M_{i,i}=1$ et $M_{j,j}=1$ alors
$x_i $ et $x_j$ sont isomorphes.  En effet, $A(x_i,x_j)$ a
$M_{ij}>0$ \'el\'ements, on peut en choisir un $f$; et
$A(x_j,x_i)$ a $M_{ji}>0$ \'el\'ements, choisissons-en $g$. Alors
$fg=1$ et $gf=1$ car $|A(x_i,x_i)|=M_{ii}=1$ et
$|A(x_j,x_j)|=M_{jj}=1$. Alors $A$ est non-r\'eduite contradiction
donc $M$
ne marche pas.\\
\begin{theoreme}[Leinster \cite{BergerLeinster}]
\label{leinsterthm} Soit $M=(m_{ij})$ une matrice carr\'ee dont
les coificients sont des entiers naturels et pour tout i
$m_{ii}\geq 2$ ,alors Cat(M)$\neq{\Oe}$ (i.e.d il existe une
cat\'egorie associ\'e \`a M).
\end{theoreme}
En effet: Soit $M=(m_{ij})$ de taille avec $m_{ii}\geq 2$, on pose
$n_{ij}:= m_{ij}$ pour $i\neq j$ et $n_{ii}:= m_{ii}-1$. On peut
d\'efinir une semicat\'egorie A associ\'e \`a $N$ dont les objets
sont $1,2,.....,n$, pour tout couple $(i,j)$ on a une fl\`eche
$\Phi_{ij}$ :i $\rightarrow$ j tel que $\Phi_{ij} \neq 1_{ii}$,la
loi de composition d\'efinit par si $f:i \rightarrow j$ et
$g:j\rightarrow k$ $\Phi_{ij}$ alors $gf=\Phi_{ik}$. Ensuite on
peut d\'efinir une cat\'egorie $B$ en rajoutant \`a $A$ les
identit\'es, pour tout $i$ on a $1_{ii}: i \rightarrow i$.  La
matrice de $B$ est $M$.
\begin{corollary}
Pour toute matrice positive on peut \'etudier si
cette matrice marche ou non.
\end{corollary}
En effet: soit $M=(m_{ij})$ de taille n alors il y a deux cas :\\
a- $m_{ii}>1$ pour tout $ i \in\{1,2,...,n\}$\\
b- il existe au moins  une i' tel que $m_{i'i'}=1$\\
Cas (a):\\
on a $Cat(M)\neq{\Oe}$ d'apr\'es le th\'eor\`eme
pr\'ec\'edant.\\
Cas(b):\\
1- s'il existe une seule i' tel que $m_{i'i'}=1$ ,on peut
\'etudier
cette matrice d'apres le th\'eor\`eme \ref{gen1}.\\
2- s'ils existent i,j,......l tel que
$m_{ii}=m_{jj}=........=m_{ll}=1$ alors il y a deux cas:\\
-si M une matrice r\'eduite alors M ne pas marche d'apr\'es (lemme \ref{marchepas}) .\\
-si M une matrice non r\'eduite alors il existe une matrice N
r\'eduction de M facile \`a \'etudier par r\'ecurrence.\\
\\
\begin{corollary}
Si $M$ est une matrice de taille $n\geq 3$ avec $m_{ij}\geq 1$,
alors $Cat(M)\neq \Oe$ si et seulement si, pour toute sous-matrice
$N\subset M$ de taille $3$ on a $Cat(N)\neq \Oe$.
\end{corollary}
En effet, dans l'\'etude pr\'ec\'edente, dans les cas (a) et (b1)
les conditions ne concernent que les triples d'indices $i,j,k$ et
donc ne concernent que les sous-matrices de taille $3$. Pour cas
(b2) si $m_{ii}=m_{jj}=........=m_{ll}=1$, la condition
n\'ecessaire et suffisante pour que $M$ marche est que pour tout
autre $k$ on a $m_{ik}= m_{jk} = ... = m_{lk}$ et $m_{ki}= m_{kj}
= ... = m_{kl}$, et que la sous-matrice d\'efinie en enlevant
$j,...,l$ marche d'apr\`es le cas (b1).

\end{document}